\documentstyle[amstex,amsfonts,amssymb,amscd,12pt]{article}

\newcommand{\QQ }{\Bbb Q}
\newcommand{\CC }{\Bbb C}
\newcommand{\ZZ }{\Bbb Z}

\newtheorem{th}{Theorem}

\newtheorem{def/th}{Definition/Theorem}

\begin{document}

\author{Bumsig Kim\thanks{Partially supported by a Sloan Doctoral
Dissertation Fellowship} }
\title{
On Equivariant Quantum Cohomology}

\date{8-2-96}
\maketitle

\section{Introduction}

There is exactly one straight line passing through any two given
distinct points; there is exactly one quadratic curve on the complex
projective plane passing through $5$ given generic points ...
One can formulate many similar enumerative problems about compact
holomorphic curves in K\"ahler manifolds, and some of such problems  have been
subject to intensive study by algebraic geometers since the last century.
Recently it was found that answers to these questions give rise
to {\em symplectic} invariants of the manifolds. First, Gromov \cite{Gr}
% Pseudo-holomorphic curves in symplectic manifolds, Inventiones Mathematikae
% (1985), ?? - ??
suggested to use pseudo-holomorphic curves in symplectic
manifolds in order to distinguish equivalence classes of symplectic structures.
Then Floer \cite{Fl}
% Symplectic fixed points and holomorphic spheres, Commun. Math. Phys. 1987
applied this idea to Arnold's fixed point conjecture \cite{Ar}
% V. Arnold, Comptes Rendus Acad. Sci. Paris, (1965)
and introduced what is called now {\em the
quantum cohomology algebra} of a symplectic manifold. Eventually, when
the ideas of symplectic topology merged with those of conformal field theory,
it became clear that various {\em Gromov -- Witten
invariants} of symplectic manifolds responsible for enumeration of holomorphic
curves can be actually calculated within the quantum cohomology algebras. Thus
numerous enumerative questions were reduced to computation of the quantum
cohomology algebras themselves.

The quantum cohomology algebra of a given compact K\"ahler manifold $X$ is,
by definition, the cohomology space of the manifold provided with a new
multiplication. Given three cohomology classes represented by cycles $a,b,c$
Poincar\'e-dual to them, one can think of the structural constants
$<a \cdot b, c>$ of the ordinary cup-product as of the number of points of
transversal intersection $a\cap b\cap c$ counted with appropriate signs.
Similarly, the structural constant $<a * b, c>$ of the quantum multiplication
can be understood as the algebraic number of isolated solutions to the
following enumerative problem:

\medskip

{\em find the number of holomorphic maps $\CC P^1 \to X$ with the points
$0,1,\infty \in \CC P^1$ mapped to the given cycles $a,b,c$, respectively.}

\medskip

Here the maps of different degrees should be counted separately, and the
degree $d$ maps contribute to the structural constant by $\pm q^d$ (so that
the ``constants'' are polynomials in $q$). Since the degree $0$
holomorphic maps are constant, the quantum multiplication turns out to be
a deformation of the ordinary cup-product: $<a*b,c>|_{q=0}=<a\cdot b, c>$.
The number of parameters of the deformation is equal to the rank $l$ of the
second homotopy group of the target K\"ahler manifold: the homotopy
class $d=:(d_1,...,d_l)\in \ZZ ^l$ of the holomorphic map is represented by
the monomial $q^d:= q_1^{d_1}...q_l^{d_l}$.

The quantum multiplication provides a machinery for answering various
enumerative problems. In particular, {\em the algebraic numbers of holomorphic
maps $(\CC P^1, z_1,...,z_N) \to (X, a_1, ..., a_N) $ of all degrees are
coefficients in the $q$-polynomial $(a_1*...*a_N, [X])$}. This follows
from the structural properties of Gromov-Witten invariants (called
sometimes {\em the composition rules}) which originate from the axioms of
Topological Field Theory (TFT) \cite{At}.

A rigorous construction of quantum cohomology algebras of general compact
symplectic manifolds including a proof of the axioms of TFT
constitutes a highly non-trivial mathematical problem
whose solution took several years (see \cite{RT, MS, Ko, BM, Beh, LT, RT2}).

\bigskip

In the present paper we study quantum cohomology algebras of flag manifolds.

\bigskip
Let $F_{s_0,...,s_l}$ be the manifold of all flags
\[ 0\subset \CC ^{s_0}\subset ...\subset \CC ^{s_{l-1}}\subset \CC ^{s_l}
=\CC ^n \]
of complex linear subspaces in $\CC ^n$ of dimensions
$0< s_0 < s_1 ... < s_l $. For an $m$-dimensional complex vector bundle
with Chern classes $c_1,...,c_m$ introduce the Chern polynomial
$x^m+c_1x^{m-1}+...+c_m$. Denote $P_0,...,P_l$ the Chern polynomials of
the tautological bundles of dimensions
$k_0=s_0, k_1=s_1-s_0,...,k_l=s_l-s_{l-1}$ over the flag manifold with the
fibers $\CC ^{s_0}, \CC ^{s_1}/\CC ^{s_0},...,\CC ^{s_l}/\CC ^{s_{l-1}}$,
respectively. The cohomology algebra $H^*(F_{s_0,...,s_l}, \QQ )$ of
the flag manifold is multiplicatively generated by the Chern classes
$(c_j^{(i)}), i=0,...,l, j=1,...,k_i$ of these bundles. A complete
set of relations between these multiplicative generators can be written
in the elegant form of a single relation between the Chern polynomials:
\[ P_0(x) ... P_l(x) = x^n.\]
We describe below the quantum deformation
of this formula.

The fraction
\[ P_0 + \frac{(-1)^{k_0+1}q_1}{P_1+
\frac{(-1)^{k_1+1}q_2}{P_2 + \frac{...}{ ... 
+\frac{(-1)^{k_{l-1}+1}q_l}{P_l}}}} \]
can be written unambiguously as the ratio $P(x)/Q(x)$ of polynomials
of degree $n$ and $n-s_0$, respectively. The coefficients of the polynomial
\[ P=x^n+\Sigma_1 x^{n-1} + ... + \Sigma_n \]
are polynomial expressions in the letters $(c_j^{(i)})$ and $q_1,...,q_l$.
The quantum deformation of the above relation reads $P=x^n$.

\bigskip

\begin{th}\label{quant-partial}
 The quantum cohomology algebra of the flag manifold
$F_{s_0,...,s_l}$ is multiplicatively generated by the $n$ Chern classes
$(c_j^{(i)})$ and the $l$ parameters $q_1,...,q_l$
satisfying the relations $\Sigma_1 =0, ...,
\Sigma_n=0$. The Poincar\'e intersection index $(a, b)$ between two
cohomology classes represented in the quantum cohomology algebra
by the polynomials $a(c, q), b(c, q)$ of these
$n+l$ generators is given by the $n$-dimensional residue
\[ <a,b> (q)=(\frac{1}{2\pi \sqrt {-1}})^n
\oint _{|\Sigma _i|=\epsilon _i} a(c,q)b(c,q)
\frac{ d c_1^{(1)}\wedge ... \wedge dc_{k_l}^{(l)}}{\Sigma_1(c,q) ...
\Sigma_n(c,q)} \ .\]
\end{th}

\bigskip
{\bf Corollary}{\it
For $N$ given generic cycles $a_1, ..., a_N $ of real
codimension $2$ in the flag manifold, the number of degree $d$
holomorphic maps
\[ (\CC P^1, x_1,...,x_N) \to (F_{s_0,...,s_l}, a_1,...,a_N) \]
is equal to the above residue with $a=a_1(c)a_2(c)...a_N(c)$ and $b=1$
where $a_{\alpha }(c)$ is the linear combination of classes $c_1^{(i)}$
Poincar\'e-dual to the cycle $a_{\alpha }$.}

\bigskip

{\em Remarks.} 1) Denote $p_i, i=1,...,l$, the $1$st Chern class
$c_1^{(i)}+...+c_1^{(l)}$ of the tautological {\em quotient} bundle
with the fiber $\CC ^n/\CC ^{s_{i-1}}$. The geometrical meaning of the
parameters $q_i$ in the above formulation is determined by the following
convention: the monomial $q_1^{d_1}...q_k^{d_k}$ represents holomorphic
curves $C$ in the flag manifold with $\int _C p_i = d_i$.

2) The statement of theorem \ref{quant-partial}
was conjectured independently by Astashkevich -- Sadov
\cite{AS} and Kim \cite{Ki}
and first proven by Ciocan-Fontanine \cite{C-F} 
for the case of complete flag manifolds
$F_{1,2,...n}$.
In the special case of manifolds $F_{1,2,...n}$ of complete flags in $\CC ^n$
it was conjectured by Givental -- Kim \cite{GK} in the form of a surprising
relation with complete integrable systems. Namely, the polynomials
$\Sigma _1(c,q),...,\Sigma_n(c,q)$ turn out to be Poisson-commuting
conservation laws of the Toda lattice (see \cite{GK} for details).
The statement of theorem \ref{quant-partial}
 in the case of Grassmannians $F_{k,n}$
is due to Witten and Siebert -- Tian \cite{W, ST} and was conjectured by Gepner.
For complex projective spaces $F_{1,n}$ the first computation
of the quantum cohomology algebra can be found in \cite{FW}.

\bigskip

The heuristic proof of theorem \ref{quant-partial}
 suggested in \cite{GK} (for complete
flags), \cite {Ki} and (with slight modification) in \cite{AS} was based on
several natural hypotheses about existence and general properties of an
equivariant generalization of the quantum cohomology theory. Given a fibration
$E\to B$ of compact manifolds with the compact K\"ahler manifold $X$ in the role
of the fiber, one can formulate various enumerative questions about holomorphic
curves in the fibers passing by marked points through given cycles in the
total space $E$ of the fibration. In particular, such a parametric
enumerative geometry can be associated, in particular,
with any principal $G$-bundle over $B$ where $G$ is a compact Lie group of
automorphisms of $X$. The enumerative information about all such bundles
can be encoded by structural constants of the {\em $G$-equivariant quantum
cohomology algebra} of $X$ which is accountable, by definition, for
enumeration of fiber-wise holomorphic spheres in the $X$-bundle $X_G\to BG$
associated with the {\em universal} principal $G$-bundle $EG\to BG$.
Like the non-equivariant quantum cohomology algebra, the equivariant one
is a deformation of the multiplicative structure in the ``classical''
equivariant cohomology $H_G^*(X):= H^*(X_G)$ in the category of algebras
over the ring $H^*(BG)$ of characteristic classes of principal $G$-bundles.

We study the $U_n$-equivariant quantum cohomology algebras of the flag
manifolds with respect to the natural action of the unitary group on $\CC^n$
and deduce theorem \ref{quant-partial}{
 from its equivariant generalization.

\bigskip

Denote $c_1,...,c_n$ the universal Chern classes of principal $U_n$ bundles.

\begin{th}\label{equant-partial}
The $U_n$-equivariant quantum cohomology
$\QQ [c_1,...,c_n]$ - algebra of the flag manifold $F_{s_0,...,s_l}$ is
isomorphic to
\[ \QQ [c_1^{(1)},...,c_{k_l}^{(l)}, q_1,...,q_l, c_1,...,c_n ]/
( \Sigma_1(c,q)-c_1, ..., \Sigma_n(c,q)-c_n )\]
The equivariant Poincar\'e pairing is given by the residue
\[ (a,b) (q,c_1,...,c_n) =(\frac{1}{2\pi \sqrt {-1}})^n
\oint _{|\Sigma_i-c_i|=\epsilon_i}
a b \frac{\wedge_{i=0}^l\wedge_{j=1}^{k_i} dc_j^{(i)}}
{\Pi _{m=1}^n (\Sigma_m (c,q)-c_m)} \ .\]
\end{th}

{\em Remark.}
1). In section 3 and 4 we construct equivariant quantum cohomology
theory for simply connected homogeneous K\"ahler spaces, 
prove the appropriate composition rule
and the other general properties of equivariant quantum cohomology
assumed in the heuristic computation in \cite{GK}, deduce theorem
\ref{equant-partial} and
obtain theorem \ref{quant-partial}
as its specialization at $c_1=0, ..., c_n=0$.

2). Constructing \lq\lq vertical
    quantum cohomology" introduced in \cite{AS},
      Lu also proved theorem \ref{equant-partial} in \cite{Lu}.

\section{Lemma}\label{eqc}
In \cite{GK, AS, Ki} the computations of quantum cohomology of
flag varieties were established with an assumption. The assumption was that
there is a $\ZZ $-graded
equivariant quantum cohomology with the properties of product,
induction, and restriction for flag varieties. Namely,

%%We prove the assumption for a generalized flag variety $X$,
%%\footnote{A generalized flag manifold is, by definition, a semisimple
%%complex Lie group quotient by a parabolic subgroup.}
%%using Kontsevich's stable maps, so that we can complete the computations
%%in \cite{GK, AS, Ki}. To do so, first we define
%%equivariant Gromov-Witten classes
%%and the equivariant quantum cohomology of $X$.
%%They should satisfy all the equivariant
%%versions of axioms in \cite{KM}, for instance, the associativity that we will
%%show. The simple idea is to make diagonally
%%equivariant all pictures in the ordinary case. And then we prove the three
%%properties noticed by Givental. For the product rule, our quantum cohomology
%%is a structure of the super commutative associative algebra with identity on
%%the tangent bundle $TH^{*}(X,\QQ )$ {\em restricted on }$H^2(X,%
%%\QQ ).$

\bigskip
{\bf Lemma} {\it
Let $G$ be a connected compact Lie group continuously acting on 
a generalized flag variety X.
Then there is a $\ZZ$-graded
equivariant quantum cohomology algebra $QH_G^{*}(X,\QQ )$ which
is $H^{*}(X_G,\QQ )\otimes _{\QQ }\QQ [q]$ as a
free $H^{*}(BG,\QQ )\otimes _{\QQ }\QQ [q]$-module, and $%
q=(q_i)$ is a formal multi-variable for a suitable
basis of $H_2(X,\ZZ )$.
The grading
is given by the usual grading on classes and the Chern number $2c_1(TX)[q_i]$
on each $q_i$. When $G$ is the trivial group, $QH_G^{*}(X,\QQ )$
  becomes the ordinary
quantum cohomology algebra.\footnote{For definition, see section 2.1.}
It has the following properties.

{\bf Product}: Let $G^{\prime }$ and $G^{\prime \prime }$ be connected
compact Lie
groups with %holomorphic
actions on $X^{\prime }$ and $X^{\prime \prime }$,
respectively. Then
\[
QH_{G^{\prime }\times G^{\prime \prime }}^{*}(X^{\prime }\times X^{\prime
\prime },\QQ )\cong QH_{G^{\prime }}^{*}(X^{\prime },\QQ )\otimes _{%
\QQ }QH_{G^{\prime \prime }}^{*}(X^{\prime \prime },\QQ ).
\]

{\bf Restriction}: Let $G^{\prime }$ be a connected
Lie subgroup of a connected compact Lie
group $G$ with a $G$-space $X$. Then, as $H^{*}(BG^{\prime },\QQ )$-algebras,
\[
QH_{G^{\prime }}^{*}(X,\QQ )\cong QH_G^{*}(X,\QQ )\otimes _{H^{*}(BG,%
\QQ )}H^{*}(BG^{\prime },\QQ ).
\]

{\bf Induction}: Let $G^{\prime }$ be a connected
Lie subgroup of a connected compact Lie
group $G$, and let $G^{\prime }$ act on $Y.$
Define $X:=G\times _{G^{\prime }}Y,$ which has the induced $G$-action$.$
Suppose $X$ becomes another generalized flag manifold with
the holomorphic quotient maps $X\rightarrow Y$ and $X\rightarrow G/G'$, then,
as $H^{*}(BG,\QQ )$-algebras,
\[
QH_G^{*}(X,\QQ )\otimes _{\QQ [q,q^{\prime }]}\QQ %
[q]\cong QH_{G^{\prime }}^{*}(Y,\QQ )
\]
where $q$ $($resp. $q^{\prime })$ is a
formal multi-variable for a suitable basis of
$H_2(Y,\ZZ )$ $($resp. $H_2(G/G^{\prime },%
\ZZ ))$. Here $q^{\prime }$ acts trivially on $q$ in $\QQ %
[q,q^{\prime }]$-module $\QQ [q]$.}

\bigskip

This lemma needs some explanations:

1. Here a suitable basis of $H_2(X,\ZZ )$ is a basis of consisting of elements 
represented by rational curves in $X$.
Let $X=G/P$, $G$ a complex Lie group (not $G$ in the lemma), $P$ a parabolic
subgroup containing a Borel subgroup $B$, $P'$ a parabolic subgroup
containing $P$ and having one more roots than $P$. Then the fibers of $G/P
\rightarrow G/P'$ are rational curves, representing an element 
of $H_2(X,\ZZ )$. Varying $P'$ provides the basis.

2. The induced map $BG'\rightarrow BG$ from the inclusion $G'\rightarrow G$
provides a natural $H^*(BG)$-module structure on $H^*(BG')$. This module
structure is used in the restriction and the induction.

3. According to a degenerating Leray spectral sequence of homology of
the fibration $X:=G\times _{G'}Y\rightarrow Y$ in the induction, 
$H_2(X,\ZZ )\cong H_2(Y,\ZZ )\oplus H_2(G/G',\ZZ )$. In this identification
the suitable basis of $H_2(X)$ decomposes into the suitable basis of $H_2(Y)$
and a basis of $H_2(G/G')$.

\bigskip

We can prove theorem \ref{quant-partial},
using this lemma, computations of
equivariant quantum cohomology of Grassmannians, and
two additional relationships: (a) the 
equivariant quantum cohomology algebra modulo $G$-characteristic classes
becomes the non-equivariant quantum cohomology, and 
(b) the equivariant quantum cohomology algebra modulo $q$'s becomes
the usual equivariant cohomology.
The proof of (a) follows from the restriction rule 
stated in the lemma,
and the proof of (b) follows from the definition of equivariant
quantum cohomology algebras given in subsection \ref{egwc}.
The computation of equivariant quantum cohomology of Grassmannians
can be obtained using Sibert-Tian's proof \cite{ST, W}, (a), and (b).
Details are in \cite{Ki}.

\bigskip

In \ref{gwc} and \ref{sa},
we collect all the facts that we are going to use to
prove the lemma. Those facts are due to
Kontsevich, Behrend and Manin, and  Pandharipande \cite{Ko, BM, Pa}.
The proof of the lemma is
presented in \ref{egwc} except for proofs of its
rules, which are in section \ref{prule}.

\section{Definition and associativity}

\subsection{Gromov-Witten classes}\label{gwc}

For a compactification of moduli space of rational maps, 
the notion of stable maps
was introduced \cite{KM, Ko}. Let $C$ be a connected, compact,
reduced, arithmetic genus zero curve $C$ with $n$ ordered marked points 
at regular points and with at most ordinary double singular points. 
%Marked points and singular points are called special points.
A stable map is a pair $(C,f)$ consisting of $C$ and a 
holomorphic map $f$ from $C$ to $X$, such that 
every irreducible component of $C$ that maps to a constant point must
have at least three special points. 
Marked points and singular points are called special points.
Let $\overline{\cal{M}}_n(X,d)$ denote the moduli
space of equivalent classes of
stable maps of degree $d\in H_2(X)$. 
Two stable maps $(C,f)$ and $(C',f')$ will be called equivalent if
there is an  
isomorphism $h$ from $C$ to $C'$ such that $f=f'\circ h$, and $h$ preserves
the ordered marked points.
The stable maps are defined to ensure that the automorphism group
of $(C,f)$ is discrete. When $X$ is a point, the moduli space becomes
the Deline-Mumford compactification $\overline{\cal{M}}_n$
of stable $n$-pointed curves of genus $0$. Note that in this
case $n$ should be greater than or equal to 3.

Let $X$ be a generalized flag variety. 
It is then shown that the moduli space  $\overline{\cal{M}}_n(X,d)$
of stable maps is
an irreducible (projective) variety with finite quotient
singularities, and
the complex dimension of the
space $\overline{\cal{M}}_n(X,d)$ is $\int_dc_1(T_X)+\dim X+n-3,$ the
``right'' dimension  \cite{Ko, BM, Pa}. 
Therefore we need not go into the difficulty of
finding a \lq virtual fundamental class'. 
According to  Behrend and Manin\cite{BM}, 
it has morphisms, a contraction $\pi ^X$, and evaluations 
at marked points: 
\begin{equation}
\begin{array}{ccc}
\overline{\cal{M}}_n(X,d) & \stackrel{ev^X}{\rightarrow } & X^n\stackrel{%
pr_i^X}{\rightarrow }X \\ 
\downarrow _{\pi ^X} &  &  \\ 
\overline{\cal{M}}_n. &  & 
\end{array}
\label{gw}
\end{equation}
After \cite{KM, Ko} let us define 
Gromov-Witten classes $%
I_{n,d}^X:H^{*}(X)^{\otimes n}\rightarrow H^{*}(\overline{\cal{M}}_n)$
in the following way: 
\[
I_{n,d}^X(a_1\otimes \cdots \otimes a_n):=(\pi ^X)_{*}(ev^X)^{*}(a_1\otimes
\cdots \otimes a_n). 
\]
In particular, this defines $I^X_{3,d},$ which 
gives a quantum multiplication structure on $%
H^{*}(X)\otimes _{\QQ }\QQ [q]$: 
there is a unique multiplication such
that 
\[
<a_1\cdot a_2,a_3>=\sum_{d\in H_2(X)}
\Pi _iq_i^{\gamma _i(d)}I^X_{3,d}(a_1,a_2,a_3),
\]
where $q=(q_1,...)$ is a formal multi-variable 
for a basis $\{\gamma _i\}$ in the closed K\"{a}hler cone, 
and $<,>$ is the $q$-linear expansion of the ordinary Poincar\'{e} pairing.
Let us choose $\{\gamma _i\}$ as the dual basis of 
the suitable basis of $H_2(X, \ZZ)$ explained in section 2.
So $\gamma _i(d)\ge 0$ for $d$ which can be represented by rational curves.
In \cite{KM}, instead of formal $q^d$, $\exp(-\int _d\omega )$ is used for a 
fixed K\"{a}hler class $\omega$. In next section we will see the
associativity of this quantum multiplication.

\subsection{The splitting axiom}\label{sa}

Let $\varphi _S:\overline{\cal{M}}_{n_1+1}\times \overline{\cal{M}}%
_{n_2+1}\rightarrow \overline{\cal{M}}_n$ be the morphism associated with
ordered partition $S=(S_1,S_2),$ $S_1\coprod S_2=\{1,...,n=n_1+n_2\},$ and 
$\varphi _S$ combines two stable curves at the 
$n_1+1$-th marked point and the first 
marked point, respectively. Let $\sum_{i,j}\eta ^{i,j}\alpha _i\otimes
\beta _j$ be the Poincar\'{e}-dual class of the diagonal $\Delta \subset
X\times X.$ The splitting axiom reads: 
\begin{eqnarray*}
&&\varphi _S^{*}(I^X_{n,d}(a_1\otimes ...\otimes a_n)) \\
&=&
\sum_{d=d_1+d_2}\sum_{i,j}I^X_{n_1+1,d_1}((\bigotimes_{k_1\in S_1}a_{k_1})%
\otimes \alpha _i)\eta ^{i,j}\otimes I^X_{n_2+1,d_2}(\beta _j\otimes
(\bigotimes_{k_2\in S_2}a_{k_2})).
\end{eqnarray*}
This axiom is proven by Behrend and Manin \cite{BM}.
In particular, the splitting axiom for $n=4$ verifies 
\begin{eqnarray*}
&&\sum_{d=d_1+d_2}\sum_{i,j}I^X_{3,d_1}(a\otimes b\otimes \alpha _i)\eta
^{i,j}I^X_{3,d_2}(\beta _j\otimes c\otimes d) \\
&=&\sum_{d=d_1+d_2}\sum_{i,j}I^X_{3,d_1}%
(b\otimes c\otimes
\alpha _i)\eta ^{i,j}I^X_{3,d_2}(\beta _j\otimes a\otimes d),
\end{eqnarray*}
which gives the associativity of quantum multiplications. Here 
one use the fact $\overline{\cal{M}}_4=\CC P^1$.

We would like to ``recall'' a proof of the splitting axiom when $n=4$:
Note that we have 
\[
\begin{array}{ccc}
\overline{\cal{M}}_{3}(X,d_1)\times \overline{\cal{M}}%
_{3}(X,d_2)\;\;\;\;\;\;\; & \overline{\cal{M}}_3(X,d) & \stackrel{ev}{%
\rightarrow }X^3 \\ 
\;\downarrow _{\pi _{d_1}\times \pi _{d_2}} & \downarrow _\pi  &  \\ 
\overline{\cal{M}}_{3}\times \overline{\cal{M}}_{3}\;\;\;%
\;\;\;\stackrel{\varphi _S}{\rightarrow } & \overline{\cal{M}}_4 & 
\end{array}
.
\]
For $d=d_1+d_2$, let $ev_{i,d_1}$ (resp. $ev_{i,d_2})$ denote
the evaluation maps from $\overline{\cal{M}}_{3}(X,d_1)$ (resp. $%
\overline{\cal{M}}_{3}(X,d_2))$ at the $i$-th marked point, where $%
i=1,2,3$. Let $\Delta $ be the diagonal in $X\times X.$
Then, from the ordered partition $S,$ 
we have the associated map $\Delta _{d_1,d_2}$
from $(ev_{3,d_1}\times ev_{1,d_2})^{-1}(\Delta )$ to $\overline{%
\cal{M}}_4(X,d),$ combining the third marked \lq point' from $\overline{%
\cal{M}}_{3}(X,d_1)$ 
with the first marked \lq point' from $\overline{%
\cal{M}}_{3}(X,d_2).$  
The variety $(ev_{3,d_1}\times ev_{1,d_2})^{-1}(\Delta )$ 
should be considered
a fibered product, and it is also an orbifold
 because $ev_i$ is a smooth morphism 
(submersion if one would like differentiable orbifold languages).
In summary, we have the following commutative
diagram of morphisms 
\[
\begin{array}{ccc}
(ev_{3,d_1}\times ev_{1,d_2})^{-1}(\Delta ) & \stackrel{\Delta
_{d_1,d_2}}{\rightarrow } & \text{Im}\Delta _{d_1,d_2} \\ 
\downarrow &  & \downarrow _\pi \\ 
\overline{\cal{M}}_{3}\times \overline{\cal{M}}_{3} & 
\stackrel{\varphi _S}{\rightarrow } & \text{Im}\varphi _S
\end{array}
. 
\]
The horizontal maps are isomorphisms because the associated trees for
stable maps are simply connected, marked by points, and labeled by degrees. 
Note that, as analytic fundamental classes, $\sum_{d_1+d_2=d}[\text{Im}\Delta
_{d_1,d_2}]=[\pi ^{-1} (\text{Im}\varphi _S)]$. 
Hence, keeping in mind that $%
\sum_{i,j}\eta ^{i,j}(ev_{3,d_1}\times ev_{1,d_2})^{*}(\alpha
_i\otimes \beta _j)$ is the Poincar\'{e}-dual class of $(ev_{3,d_1}\times
ev_{1,d_2})^{-1}(\Delta )$ in 
$\overline{\cal{M}}_{3}(X,d_1)\times 
\overline{\cal{M}}_{3}(X,d_2),$ we conclude the proof.

\subsection{Equivariant Gromov-Witten classes}\label{egwc}

Let $X$ have a 
continuous $G$-action, $G$ being connected and compact. Then we have
maps 
\[
\begin{array}{ccc}
\overline{\cal{M}}_n(X,d)\times _GEG & \stackrel{ev^{X_G}}{\rightarrow } & 
X^n\times _GEG \\ 
\downarrow _{\pi ^{X_G}} &  &  \\ 
\overline{\cal{M}}_n\times BG &  & 
\end{array}
, 
\]
the equivariant version of the diagram (\ref{gw}).
Recall that $H_{G}^{*}(X^n)$ is 
$H_G^{*}(X)^{\otimes n}$, due to the
projections $(X^n)_G\rightarrow X_G$, so that they can be identified. 
Define equivariant Gromov-Witten
classes $I_{n,d}^{X_G}:H_G^{*}(X)^{\otimes n}\rightarrow H^{*}
(\overline{\cal{M}}_n)\otimes _{\QQ }H^{*}(BG)$ by 
\[
I_{n,d}^{X_G}(a_1\otimes \cdots \otimes a_n):=\pi _*^{X_G}%
(ev^{X_G})^*(a_1\otimes
\cdots \otimes a_n), 
\]
where $a_i\in H_G^{*}(X)$.

The module $H_G^{*}(X)\otimes _{\QQ }\QQ [q]$ has a unique
multiplication by the characterization 
\[
<a_1\cdot a_2,a_3>=\sum_dq^dI^{X_G}_{3,d}(a_1,a_2,a_3),
\label{emuti}
\]
where $<,>$ is the $q$-linear expansion of the 
{\em equivariant} Poincar\'e pairing.
For the {\em equivariant} Poincar\'{e}-dual class 
$\sum_{i,j}\eta ^{i,j}\alpha _i\otimes \beta _j$
of the diagonal $\Delta _G\subset (X\times X)_G$
the equivariant version of the splitting axiom holds, namely,
\begin{eqnarray*}
&&\varphi _S^{*}(I_{4,d}^{X_G}(a_1\otimes ...\otimes a_n)) \\
&=&\sum_{d=d_1+d_2}\sum_{i,j}I_{3,d_1}^{X_G}((\bigotimes_{k_1\in
S_1}a_{k_1})\otimes \alpha _i)\eta ^{i,j}\otimes I_{3,d_2}^{X_G}
(\beta _j\otimes
(\bigotimes_{k_2\in S_2}a_{k_2})),
\end{eqnarray*}
where all tensor products are from $H^{*}(BG)$-module structures.
The point is that all maps in the proof of the
splitting axiom for the nonequivariant version are (diagonally) equivariant.
Just as in the nonequivariant case, let us keep in mind that
we have
$\Delta _{d_1,d_2}:(ev_{3,d_1}\times ev_{1,d_2})^{-1}(\Delta
)\rightarrow \overline{\cal{M}}_n(X,d)$, its equivariant version,
and that
$\sum_{i,j}\eta ^{i,j}(ev_{3,d_1}\times ev_{1,d_2})^{*}(\alpha
_i\otimes \beta _j)$ is the equivariant Poincar\'{e}-dual class of $\left(
(ev_{3,d_1}\times ev_{1,d_2})^{-1}(\Delta )\right) _G$ in $\left( 
\overline{\cal{M}}_{3}(X,d_1)\times \overline{\cal{M}}%
_{3}(X,d_2)\right) _G.$
Then the proof follows from the parallel argument of the proof of the ordinary
splitting property given in the previous section.

\bigskip
{\bf Definition/Theorem} {\it
Analogous to the one
defining the quantum cohomology, we define the
equivariant quantum cohomology multiplication $QH_G(X)$.
The associativity can be
proven by the equivariant version of the splitting property when $n=4$.
The ring is graded as stated in the lemma. }

\bigskip

When $G$ is the trivial group, obviously the equivariant quantum cohomology
is the ordinary quantum cohomology.
Since $H^*(X)$ and $H^*(BG)$ are generated by even degree classes, 
$H^*_G(X)=H^*(X)\otimes _{\QQ }H^*(BG)$ as linear spaces, and 
$QH^*_G(X)$ is a free $H^*(BG)\otimes _{\QQ } \QQ [q]$-module.

\section{Rules}\label{prule}

{\bf A proof of product rule}:

Suppose $G^{\prime }$ and $G^{\prime \prime }$ are 
connected compact Lie groups. Let $%
X^{\prime }$ be a $G^{\prime }$ space, and let $X^{\prime \prime }$ be a $%
G^{\prime \prime }$ space, then we have the induced $G^{\prime }\times
G^{\prime \prime }$ space, $X^{\prime }\times X^{\prime \prime }$ and,
as $H^*(BG^{\prime }\times
BG^{\prime \prime })=H^*(BG^{\prime })\otimes _{\QQ }H^*(BG^{\prime
\prime })$-modules, $%
H_{G^{\prime }\times G^{\prime \prime }}^{*}(X^{\prime }\times X^{\prime
\prime })\cong H_{G^{\prime }}^{*}(X^{\prime })\bigotimes_{\QQ %
}H_{G^{\prime \prime }}^{*}(X^{\prime \prime })$.
Since the complement to the subset $\overline{%
\cal{M}}_3^0(X,d)$ consisting of ``smooth'' curve is a divisor (with
normal crossings), $
\overline{\cal{M}}_3(X^{\prime },d_1)\times \overline{\cal{M}}%
_3(X^{\prime \prime },d_2)$ and $\overline{\cal{M}}_3(X^{\prime }\times
X^{\prime \prime },d)$ are birational, so that 
$I_{3,d_1}^{X^{\prime }}\cdot I_{3,d_2}^{X^{\prime \prime
}}=I_{3,(d_1,d_2)}^{X^{\prime }\times X^{\prime \prime }}.$
Let $C$ and $D$ be finite cycles of $BG^{\prime }$ and $%
BG^{\prime \prime }$ respectively. Then integrating fibers
over $C\times D$, we conclude
$I_{3,d_1}^{X_{G^{\prime }}^{\prime }}\bigotimes_{\QQ %
}I_{3,d_2}^{X_{G^{\prime \prime }}^{\prime \prime
}}=I_{3,(d_1,d_2)}^{X'\times X_{G^{\prime }\times G^{\prime \prime }}^{\prime
\prime }}.$
Hence we have the proof of the product property stated 
in the theorem.

{\bf A proof of restriction rule:}

Let $G^{\prime }\subset G$ be a Lie subgroup and $X$ be a $G$-space.
Consider $X$ a $G^{\prime }$-space for $X_G$. 
Let $p :BG^{\prime }\rightarrow BG$
be the map induced from the inclusion $G^{\prime }\subset G.$ 
We have natural induced morphisms and a diagram
\[
\begin{array}{ccc}
H^{*}(\overline{\cal{M}}_n(X,d)_{G^{\prime }}) & \leftarrow  & H^{*}(%
\overline{\cal{M}}_n(X,d)_G) \\ 
\downarrow  &  & \downarrow  \\ 
H^{*}(BG^{\prime }) &\stackrel{p^*}{\leftarrow } & H^{*}(BG).
\end{array}
\]

The diagram is commutative, since
for any finite cycle $C$ in $BG'$,
$\overline{\cal{M}}_n(X,d)\times _Gp(C)$ 
induces $\overline{\cal{M}}_n(X,d)\times _{G'}C$ by the map $p$. 
The restriction rule follows.

{\bf A proof induction rule:}

Let $G^{\prime }\subset G$ be a Lie subgroup$.$ For induction consider a
generalized flag manifold $Y$ with a $G^{\prime }$-action and let $%
X=G\times _{G^{\prime }}Y$. 
For $d\in H_2(Y)\subset H_2(X)$, there are natural identifications 
$Y_{G'}=X_G$ and
$\overline{\cal{M}}_n(Y,d)_{G^{\prime }} =\overline{\cal{M}}%
_n\left( X,d\right) _G  $. From the commutative diagram
\[
\begin{array}{ccc}
\overline{\cal{M}}_n(Y,d)_{G^{\prime }} & = & \overline{\cal{M}}%
_n\left( X,d\right) _G \\
\downarrow  &  & \downarrow  \\
BG^{\prime } &\stackrel{p}{\rightarrow}  & BG,
\end{array}
\]
$I_{n,d}^{X_G}=p_{*}I_{n,d}^{Y_{G^{\prime }}}$ and
the induction rule follows.

{\it Acknowledgement.} 
I would like to thank A. Givental for valuable suggestions, and
H. Chang for a discussion on the splitting axiom. This paper is
a part of my thesis.

Department of Mathematics, University of California, Berkeley, CA 94720, USA

{\it E-mail:} bumsig@@math.berkeley.edu

\end{document}